\documentclass[twoside, 11pt, reqno]{amsart}

\usepackage[left=2.5cm, right=2.5cm, top=3cm, bottom=2.5cm]{geometry}

\usepackage[colorlinks=true, pdfstartview=FitV, linkcolor=blue,citecolor=blue, urlcolor=blue]{hyperref}

\usepackage[usenames]{xcolor}
\definecolor{labelkey}{rgb}{0,0,1}
\definecolor{Red}{rgb}{0.7,0,0.1}
\definecolor{Green}{rgb}{0,0.7,0}

\usepackage{amsfonts, amssymb, amsmath, amsthm, mathrsfs, bbm, cjhebrew, gensymb, textcomp, mathtools, dsfont,tikz}

\usepackage[normalem]{ulem}

\usepackage{commath, setspace, subcaption, parcolumns, multirow, multicol, accents, comment, marginnote, verbatim, empheq, enumerate, stackrel, enumitem}

\usepackage{cleveref}
\usepackage{graphicx}
\usepackage{epsfig}
\usepackage{psfrag}
\usepackage{float}
\usepackage[active]{srcltx}
\usepackage[pagewise, mathlines]{lineno}

\newtheorem{Theorem}{Theorem}[section]

\newtheorem{definition}{Definition}[section]

\newtheorem{Thm}{Theorem}[section]
\newtheorem{Prop}{Proposition}[section]
\newtheorem{Lem}[Theorem]{Lemma}
\newtheorem{Rmk}{Remark}[section]

\numberwithin{equation}{section}

\newcommand\lref[1]{lemma~\ref{#1}}
\newcommand{\al}{\alpha}
\newcommand{\be}{\beta}

\newcommand{\veps}{\varepsilon}

\newcommand{\Lam}{\Lambda}
\newcommand{\si}{\sigma}

\newcommand{\tht}{\theta}

\newcommand{\RR}{\mathbb{R}}

\newcommand{\Sob}[2]{\lVert#1\rVert_{#2}}

\newcommand{\bdy}{\partial}

\newcommand{\Hdot}{\dot{H}}

\title[Inviscid and fast rotation limit for dispersive generalized SQG equations]{Inviscid and fast rotation limit for dispersive generalized SQG equations}
\author{Anuj Kumar$^{1,\dagger}$, Sweety$^{2}$, Mohit Yadav$^{3}$}
\address{$^1$Department of Mathematics,
Indian Institute of Technology Jodhpur, India}
\address{$\dagger$ corresponding author}
\email[A Kumar]{akumar241@outlook.com}
\address{$^2$ Department of Mathematics,
Malviya National Institute of Technology Jaipur, India}
\email[Sweety]{2025pma5440@mnit.ac.in}
\address{$^3$ Department of Mathematics,
National Institute of Technology Warangal, India}
\email[M Yadav]{mohityadav22063@gmail.com}
\date{\today}
\thanks{}
\begin{document}
\begin{abstract}
	We consider the dissipative generalized surface quasi-geostrophic equations with a dispersive forcing term and study the behavior of solutions in the limit of vanishing viscosity. Using the approach of relative energy inequality, we prove convergence results in different dispersive regimes that extend the results obtained by Kosloff et al. (Appl. Math. Lett. (2019) 88 243-249) to more singular active scalar equations.
\end{abstract}
\maketitle
{\noindent \small {\it {\bf Keywords: Dispersive generalized SQG equations, Weak solution, Relative energy inequality, Inviscid limit, Large dispersive forcing}
}  }

\section{Introduction}
In this paper, we consider the dissipative generalized surface quasi-geostrophic equation (gSQG) with an additional dispersive forcing term, given by
\begin{align}\label{gSQG}
    \begin{cases}
        &\bdy_t \tht+u\cdot \nabla \tht+\nu (-\Delta)^{\al}\tht+\frac{1}{\veps}\mathcal{R}_1 \tht=0,\\
       &u=\mathcal{R}^{\perp}\Lam^{\be-1}\tht=\left(\bdy_{x_2}\Lam^{\beta-2}\tht,-\bdy_{x_1}\Lam^{\beta-2}\tht\right),\\
        &\tht(x,0)=\tht_0(x),\quad x\in \mathbb{R}^2,\,t>0,
    \end{cases}
\end{align}
where $\tht(x,t)$ denotes the evolving scalar, such as the surface temperature of the ocean, $u(x,t)$ denotes the advecting velocity field, and $\Lam$ denotes the fractional Laplacian $(-\Delta)^{1/2}$. $\mathcal{R}$ represents the usual Riesz transform.  $\al \in (0,1)$ denotes the dissipation parameter and $\beta\in (0,2)$ denotes the constitutive law parameter. The domain is assumed to be $\mathbb{R}^2$. The dispersive forcing term $\mathcal{R}_1 \tht=\bdy_{x_1}\Lam^{-1}\tht$ represents the advection of a large-scale buoyancy resulting from the Coriolis force \cite{Lapeyre2017}. The parameter $\varepsilon>0$ is the non-dimensional Rossby number, which measures the weight of the Coriolis force in geophysical models.
\par
In the absence of dispersive forcing, \eqref{gSQG} was first studied in \cite{ChaeConstantinWu2012} and \cite{ChaeConstantinCordobaGancedoWu2012}, while the inviscid counterpart was studied in \cite{ChaeConstantinWu2011}. For $\be\in [0,1]$, the equations interpolate between the 2D Euler equation in vorticity form ($\be=0$) and the dissipative surface quasi-geostrophic equation (SQG) ($\be=1$). The regime defined by $\be \in (1,2)$, considered in this paper, defines active scalar equations with even more singular constitutive laws than the SQG equation. The SQG equation ($\be=1$)  has attracted a lot of interest, especially due to structural analogies to the 3D Euler equation (when $\nu=0$) and the 3D Navier-Stokes equations (when $\al=1/2$). For more on SQG-related works, we refer the reader to \cite{ConstantinTarfuleaVicol2015} and the references therein. As a model for more singular constitutive laws, the gSQG equations have also been studied intensively in recent periods, and we refer the reader to \cite{ChaeConstantinCordobaGancedoWu2012, HuKukavicaZiane2015, JollyKumarMartinez2020a} and the references therein for important works on gSQG equations. 
\par
The dispersive SQG equation is an important model for studying wave-turbulence interactions as observed in various numerical studies \cite{SukhatmeSmith2009} as well as due to its analogies to the Navier-Stokes-Coriolis (NSC) system. The NSC system is a basic geophysical model describing large-scale phenomena \cite{CheminDesjardinsGallagherGrenier2006book}. Unlike the Navier-Stokes equations, the presence of the Coriolis forcing in the NSC system affords a strong stabilizing effect in the regime of fast rotation. As a result, global well-posedness is known for the 3D NSC system even for large initial data \cite{CheminDesjardinsGallagherGrenier2006book}. More specifically, Chemin et al. \cite{CheminDesjardinsGallagherGrenier2006book} employ Strichartz estimates and establish that the dispersive term stabilizes the NSC system toward a 2D Navier-Stokes type system.
The study of dispersive SQG was initiated by Kiselev and Nazarov \cite{KiselevNazarov2010}, wherein they established the existence of global smooth solutions in the torus for the case of critical dissipation ($\al=1/2$). Cannone et al. established analogous results to those of \cite{CheminDesjardinsGallagherGrenier2006book} for the supercritical dispersive SQG equation. In subsequent works of \cite{ElgindiWidmayer2015} and \cite{WanChen2016}, the dispersive and Strichartz estimates were further improved to show global well-posedness of strong solutions of the inviscid dispersive SQG equation. We crucially use the estimates developed in \cite{ElgindiWidmayer2015} (see \cref{P:1}).
\par
Our goal in this paper is to analyze the behavior of solutions to \eqref{gSQG} for $\be \in (1,2)$ in the limit of vanishing viscosity and large dispersive forcing, i.e., when $\nu$ and $\veps$ tend to 0. In the case of the dispersive SQG equation, this study was carried out by Kosloff et al. \cite{KosloffNichePlanas2019}. The main tool they used was the relative energy inequality, which was motivated by its earlier applications in studying the inviscid incompressible limit of the NSC system \cite{FeireislJinNovotny2012, CaggioNecasova2017}. This method also allowed them to reproduce the previous results describing the vanishing viscosity limit for the unforced SQG \cite{Berselli2002, Wu1997a}. To estimate the nonlinear terms, they crucially exploit the $L^p$ boundedness of the Riesz transform ($1\le p<\infty$) in the constitutive law. For gSQG equations, due to more singular constitutive laws, such an approach is no longer feasible, and a more delicate analysis by using commutators is required. Our approach is also based on the tool of the relative energy inequality as done in \cite{KosloffNichePlanas2019}; however, for dealing with the nonlinear terms, it involves a delicate approach of decomposing the terms into terms involving commutators as done previously in several works on gSQG (see \cite{ChaeConstantinCordobaGancedoWu2012, HuKukavicaZiane2015,   JollyKumarMartinez2020a, JollyKumarMartinez2020b}). We then invoke commutator estimates in \cite{Li2019} (see \lref{L:1}) to bound these terms.
To the best of our knowledge, this is the first work on the vanishing viscosity limit for the most singular regime of the dispersive gSQG equations. 

\section{Main results}
\subsection{Inviscid and fast rotation limit}
We now state our first result, which describes the behavior of \eqref{gSQG} in the limit of vanishing viscosity and large dispersive forcing. The dispersive part of \eqref{gSQG} is given by
\begin{align}\label{gSQG:dispersive:part}
    \begin{cases}
        &\bdy_t \phi+\frac{1}{\veps}\mathcal{R}_1 \phi=0,\\
        &\phi(x,0)=\phi_0 (x),\,x\in \mathbb{R}^2,\,t>0,
    \end{cases}
\end{align}
Our main result in this setting is stated as follows.
\begin{Thm}\label{T:1}
    Let $\be \in (1,2)$. Let $m>1+\be$ and $\{\tht_{0}^{\nu,\veps}\} \in L^2(\mathbb{R}^2)$ and $\{\phi_{0}^{\veps}\}\in H^m(\mathbb{R}^2)$ be uniformly bounded in $\veps$ and $\nu$ such that $\nabla \phi_{0}^{\veps}\in \dot{B}^{1+\be}_{1,1}(\mathbb{R}^2)$. Let $\tht^{\nu, \veps}$ be a weak solution to \eqref{gSQG} and $\phi^{\veps}$ be a smooth solution to \eqref{gSQG:dispersive:part}. Suppose that 
    \[\Sob{\tht_{0}^{\nu,\veps}-\phi_{0}^{\veps}}{L^2(\mathbb{R}^2)}\rightarrow 0\quad \text{as}\quad \nu,\veps \rightarrow 0.\]
    Then, we have
    \[\Sob{\tht^{\nu,\veps}-\phi^{\veps}}{L^{\infty}(0,T;L^2)\cap L^2(0,T;\Hdot^{\al})}\rightarrow 0 \quad \text{as}\quad \nu,\veps \rightarrow 0.\]
\end{Thm}
\subsection{Inviscid limit with constant dispersive forcing}
The second result we obtain describes the behavior of \eqref{gSQG} in the limit of vanishing viscosity but with constant dispersive forcing. Consider the following limiting system:
\begin{align}\label{gSQG:inviscid:part}
    \begin{cases}
        &\bdy_t \bar{\tht}+\bar{u}\cdot \nabla \bar{\tht}+\nu (-\Delta)^{\al}\bar{\tht}+\frac{1}{\veps}\mathcal{R}_1 \bar{\tht}=0,\\
       &\bar{u}=\mathcal{R}^{\perp}\Lam^{\be-1}\bar{\tht},\\
        &\bar{\tht}(x,0)=\bar{\tht}_0(x),\quad x\in \mathbb{R}^2,\,t>0,
    \end{cases}
\end{align}
The dispersive term does not contribute to energy decay (see \cref{R:3}). As a result, for initial data in $H^{m}(\mathbb{R}^2)$ with $m>1+\be$, the existence of local smooth solutions $\bar{\tht}\in C^{1}(0,T^*;H^m)$ can be obtained by a simple modification of the arguments given in \cite{ChaeConstantinCordobaGancedoWu2012, HuKukavicaZiane2015} for the non-dispersive inviscid gSQG equations where $T^*$ is the maximal local existence time. Our main result in this setting is stated as follows.
\begin{Thm}\label{T:2}
    Let $\be \in (1,2)$. Let $\veps>0$ be fixed. Let  $m>1+\be$ and $\bar{\tht}_0 \in H^m(\mathbb{R}^2)$ and $\tht^\nu_0\in L^2(\mathbb{R}^2)$. Let $\bar{\tht}$ be a local smooth solution to \eqref{gSQG:inviscid:part} and $\tht^\nu$ a global weak solution to \eqref{gSQG}. Suppose that
    \[\Sob{\tht^{\nu}_0-\bar{\tht}_0}{L^2(\mathbb{R}^2)}\rightarrow 0,\quad \text{as}\quad \nu\rightarrow 0.\]
    Then, we have
    \[\Sob{\tht^\nu-\bar{\tht}}{L^\infty (0,T^*;L^2)}\rightarrow 0, \quad \text{as}\quad \nu\rightarrow 0. \]
\end{Thm}
The proof of \cref{T:2} is similar to that of theorem 2 in \cite{KosloffNichePlanas2019}. In particular, the approach does not require dispersive estimates to control the nonlinear term. This illustrates that a fixed dispersive forcing does not contribute to energy decay and hence does not affect the inviscid limit.
\subsection{Inviscid limit with no dispersion}
The third result we obtain treats the case with no dispersive forcing ($\veps=\infty$), i.e., we consider 
\begin{align}\label{gSQG:non:dispersive}
    \begin{cases}
        &\bdy_t \tht+u\cdot \nabla \tht+\nu (-\Delta)^{\al}\tht=0,\\
       &u=\mathcal{R}^{\perp}\Lam^{\be-1}\tht,\\
        &\tht(x,0)=\tht_0(x),\quad x\in \mathbb{R}^2,\,t>0,
    \end{cases}
\end{align}
We study the relation between solutions to the dissipative gSQG equations \eqref{gSQG:non:dispersive} and their inviscid counterpart, which we denote by 
$\tht^\nu$ and $\tilde{\tht}$ respectively. For initial data in $L^2(\mathbb{R}^2)$, it can be shown that \eqref{gSQG:non:dispersive} has global weak solutions satisfying the energy inequality (see \cref{R:2}).
 Our main result in this setting is stated as follows.
\begin{Thm}\label{T:3}
    Let $\be \in (1,2)$. Let  $m>1+\be$ and $\tilde{\tht}_0 \in H^m(\mathbb{R}^2)$ and $\tht^\nu_0\in L^2(\mathbb{R}^2)$. Let $\tht^\nu$ be a global weak solution to \eqref{gSQG} and $\tilde{\tht}$ be a local smooth solution to the corresponding inviscid problem. Suppose that
    \[\Sob{\tht^{\nu}_0-\tilde{\tht}_0}{L^2(\mathbb{R}^2)}\rightarrow 0,\quad \text{as}\quad \nu\rightarrow 0.\]
    Then, we have
    \[\Sob{\tht^\nu-\tilde{\tht}}{L^\infty (0,T^*;L^2)}\rightarrow 0, \quad \text{as}\quad \nu\rightarrow 0. \]
\end{Thm}
\section{Mathematical Preliminaries}
In this section, we present the relevant notations and definitions used in the paper. For more details, the reader can refer to \cite{BahouriCheminDanchinBook2011}.
Let $\mathscr{S}(\RR^2)$ be the space of Schwartz class functions on $\RR^2$ and $\mathscr{S}'(\RR^2)$ denote the space of tempered distributions. We denote by $\widehat{f}$ or $\mathcal{F}(f)$, the Fourier transform of $f$, which is defined by
	\[\widehat{f}(\xi)\overset{}{:=}\int e^{-2\pi i x\cdot \xi}f(x)dx,\quad f\in\mathscr{S}'(\RR^2).\] 
  Let $\Lam^\si$ denotes the fractional Laplacian operator which is defined as
	    \begin{align}\notag
	        \mathcal{F}(\Lam^\si f)(\xi)=|\xi|^\si\mathcal{F}(f),\quad \si\in\RR.
	    \end{align} 
   The $L^2$-based Sobolev spaces are defined as follows. For $\si \in \mathbb{R}$, we have
	    \begin{align}
	        &\Hdot^\si(\mathbb{R}^2){:=}\left\{f\in \mathscr{S}'(\RR^2):\widehat{f}\in L^2_{loc},\, \Sob{f}{\Hdot^\si}:=\Sob{\Lam^\si f}{L^2}<\infty\right\},\label{def:hom:Sob:norm}\\
	        &H^\si(\mathbb{R}^2){:=}\left\{f\in \mathscr{S}'(\RR^2):\widehat{f}\in L^2_{loc}, \Sob{f}{H^\si}{=}\Sob{(I-\Delta)^{\si/2}) f}{L^2}<\infty\right\}.\label{def:inhom:Sob:norm}
	    \end{align}
\par
We recall the definition of a weak solution of \eqref{gSQG}:
\begin{definition}\label{D:1}
A function $\tht$ is said to be a global weak solution to \eqref{gSQG} with the initial condition $\tht_0\in L^2(\mathbb{R}^2)$ if for all $T>0$, we have
\[\tht \in L^\infty(0,T;L^2)\cap L^2(0,T;\dot{H}^\alpha)\]
and
\[\int_{\mathbb{R}^2}\tht \phi\,dx\bigg{|}_0^t=\int_0^t \int_{\mathbb{R}^2}(\bdy_s \phi)\tht+(u\cdot \nabla \phi)\tht-\nu \Lam^{\al}\phi \Lam^{\al}\tht-\frac{1}{\veps}(\mathcal{R}_1 \tht)\phi\,dx\,ds\]
for all test functions $\phi \in C_c^\infty ([0,\infty)\times \mathbb{R}^2)$ and almost every $t>0$. 
\end{definition}

\begin{Rmk}\label{R:2}
    Using the classical Leray's proof for the weak $L^2$ solutions of the Navier-Stokes equations modified for the gSQG equations analogous to the approach in \cite{ChaeConstantinCordobaGancedoWu2012}, it can be shown that there exist global weak solutions (in the sense of distributions) as in \cref{D:1}, which satisfy the energy inequality uniformly in $\veps>0$,
    \begin{align}\label{energy:inequality}
    \Sob{\tht(t)}{L^2}^2+2\nu \int_0^t \Sob{\Lam^{\al}\tht(s)}{L^2}^2\,ds\le \Sob{\tht(0)}{L^2}^2\quad \forall \,t>0.
    \end{align}
\end{Rmk}
\begin{Rmk}\label{R:3}
    Observe that the dispersive forcing term in \eqref{gSQG} does not contribute to energy decay, since as $|\widehat{\tht}(\xi)|^2=\widehat{\tht}(\xi)\widehat{\tht}(-\xi),$ we have by Plancherel's formula
    \[\langle \Lam^{\sigma}\mathcal{R}_1 \tht, \Lam^\sigma \tht\rangle_{\mathbb{R}^2}=\int_{\mathbb{R}^2}\Lam^{\sigma}\mathcal{R}_1 \tht(x)\Lam^{\sigma}\tht(x)\,dx=-\int_{\mathbb{R}^2}i\xi_1 |\xi|^{2\sigma-1}|\widehat{\tht}(\xi)|^2\,d\xi=0.\]
\end{Rmk}
We recall the following decay estimate from \cite{ElgindiWidmayer2015} for solutions to \eqref{gSQG:dispersive:part}.
\begin{Prop}[\cite{ElgindiWidmayer2015}]\label{P:1} For all $t\in \mathbb{R}$ and $f\in \dot{B}^{2}_{1,1}$, it holds that
\[\Sob{e^{\frac{1}{\veps}\mathcal{R}_1 t}f}{L^\infty}\le C \left( \frac{\veps}{\veps+|t|}\right)^{1/2}\Sob{f}{\dot{B}^{2}_{1,1}}.\]
\end{Prop}
As remarked earlier, to control the nonlinear terms, we need to decompose them into components involving commutators. To this end, we recall a commutator estimate from \cite{Li2019} which we invoke in the following section. The commutator of two operators, $S$ and $T$, is denoted by $[S,T]$, where
    \[[S,T]:=ST-TS.\]
For $s>0$, let $A^s$ denote a differential operator such that its symbol $\widehat{A^s}(\xi)$ is a homogeneous function of degree $s$ and $\widehat{A^s}(\xi)\in C^{\infty}(\mathbb{S}^{d-1})$. Then, the following result holds. In particular, it holds for operators satisfying $\widehat{A^s}(\xi)=i|\xi|^{s-1}\xi_1$.
\begin{Lem}[\cite{Li2019}]\label{L:1} For $s\in (0,1)$, $p \in (1,\infty)$, and $f,g \in \mathscr{S}(\RR^d)$, we have
\begin{align}
    \Sob{[A^s,g]f}{L^p}\le C\Sob{f}{L^p}\Sob{\Lam^s g}{BMO}\le C\Sob{f}{L^p}\Sob{\Lam^s g}{L^\infty}.
\end{align}
\end{Lem}
\subsection{Relative energy inequality}
To prove our main results, we employ the method of relative energy inequality as done in \cite{KosloffNichePlanas2019}. We start by defining the relative energy functional:
\begin{definition}\label{D:2}
    Let $\tht$ be a global weak solution to \eqref{gSQG}. Let $\phi\in C_{c}^\infty ([0,T]\times \mathbb{R}^2)$. The relative energy functional is defined by
    \[\mathcal{E}(\tht,\phi):=\int_{\mathbb{R}^2}\frac{1}{2}|\tht-\phi|^2\,dx.\]
\end{definition}
We require the following inequality, which was established in \cite{KosloffNichePlanas2019}. For the convenience of the reader, we provide the proof in the appendix.
\begin{Prop}[\cite{KosloffNichePlanas2019}]\label{P:2}
    For all $t\in [0,T]$, we have
    \[\mathcal{E}(\tht,\phi)\bigg{|}_0^t+\nu \int_0^t\int_{\mathbb{R}^2}|\Lam^{\al}(\tht-\phi)|^2\,dx\,ds\le \int_0^t {Rem}(\tht,\phi)\,ds,\]
where
\[{Rem}(\tht,\phi):=\int_{\mathbb{R}^2}(\bdy_s \phi+u\cdot \nabla \phi)(\phi-\tht)+\nu \Lam^{\al}\phi(\Lam^\al \phi-\Lam^\al \tht)+\frac{1}{\veps}(\mathcal{R}_1 \tht)(\phi-\tht)\,dx.\]
\end{Prop}

\section{Proof of Theorems}
\subsection{Proof of \cref{T:1}}
For ease of notation, we drop the superscripts $\nu$ and $\veps$. Let $\tht$ denote a weak solution to \eqref{gSQG} and $\phi$ denote a solution to \eqref{gSQG:dispersive:part}. We want to obtain an estimate for
\begin{align}\label{Remainder}
    \int_0^t Rem(\tht,\phi)\,ds=\int_0^t \int_{\mathbb{R}^2}(\bdy_s \phi+u\cdot \nabla \phi)(\phi-\tht)+\nu \Lam^{\al}\phi(\Lam^\al \phi-\Lam^\al \tht)+\frac{1}{\veps}(\mathcal{R}_1 \tht)(\phi-\tht)\,dx\,ds.
\end{align}
We estimate the second and the third term on the right just like in \cite{KosloffNichePlanas2019} to obtain
\begin{align}\label{dissipative:term:bound}
    \nu \int_0^t \int_{\mathbb{R}^2}\Lam^{\al}\phi(\Lam^\al \phi-\Lam^\al \tht)\,dx\,ds=\nu\int_0^t\int_{\mathbb{R}^2}\Lam^{2\al}\phi(\phi-\tht)\,dx\,ds&\le \nu \int_0^t\Sob{\Lam^{2\al}\phi}{L^2}\Sob{\phi-\tht}{L^2}\,ds \notag\\
    &\le \nu C\int_0^t \Sob{\phi-\tht}{L^2}\,ds,
\end{align}
where we used the fact that since $0<\al<1$ and $m>1+\be$, we have
\[\Sob{\Lam^{2\al} \phi}{L^2}\le \Sob{\phi}{H^{2\al}}\le \Sob{\phi}{H^m}\le C\Sob{\phi_0}{H^m}\le C.\]
Also, we have
\begin{align}
    &\int_0^t \int_{\mathbb{R}^2}\left(\bdy_s \phi+\frac{1}{\veps}(\mathcal{R}_1 \tht)\right)(\phi-\tht)\,dx\,ds\notag\\&=\int_0^t\int_{\mathbb{R}^2}\left(\bdy_s\phi+\frac{1}{\veps}(\mathcal{R}_1 \phi)\right)(\phi-\tht)-\frac{1}{\veps}(\mathcal{R}_1 (\phi-\tht))(\phi-\tht)\,dx\,ds\notag\\&=0,
\end{align}
where the first term on the right is 0 since $\phi$ solves \eqref{gSQG:dispersive:part} and the second term is 0 by \cref{R:3}. For the remaining term, we observe that
\begin{align*}
    \int_{\mathbb{R}^2}(u\cdot \nabla \phi)(\phi-\tht)\,dx&=\int_{\mathbb{R}^2}\nabla^{\perp}\tht\cdot \nabla \phi (\phi-\tht)\,dx \notag \\ &=\int_{\mathbb{R}^2}\nabla^{\perp}\Lam^{\be-2}\phi \cdot \nabla \phi (\phi-\tht)\,dx-\int_{\mathbb{R}^2}\nabla^{\perp}\Lam^{\be-2}(\phi-\tht)\cdot \nabla \phi (\phi-\tht)\,dx\\
    &=I+J.
\end{align*}
We have
\begin{align}\label{eq:commutator:decom:1}
    I=&\int_{\mathbb{R}^2}\nabla^{\perp}\Lam^{\be-2}\phi \cdot \nabla \phi (\phi-\tht)\,dx-\int_{\mathbb{R}^2}\Lam^{\be-2}\underbrace{(\nabla^{\perp}\phi \cdot \nabla \phi)}_{=0} (\phi-\tht)\,dx\notag\\
    =&-\int_{\mathbb{R}^2}\left([\nabla^{\perp}\Lam^{\be-2}\cdot, \nabla \phi]\phi\right)(\phi-\tht)\,dx
\end{align}
and 
\begin{align}\label{eq:commutator:decom:2}
    J&=-\int_{\mathbb{R}^2} \nabla^{\perp}\Lam^{\be-2}(\phi-\tht)\cdot \nabla {\phi}({\phi}-\tht)\,dx\notag\\
    &=\int_{\mathbb{R}^2} ({\phi}-\tht)(\nabla^{\perp}\Lam^{\be-2}\cdot( \nabla {\phi}({\phi}-\tht)))\,dx \notag\\
    &=\frac{1}{2}\int_{\mathbb{R}^2}\left([\nabla^{\perp}\Lam^{\be-2}\cdot,\nabla {\phi}]({\phi}-\tht)\right)({\phi}-\tht)\,dx.
\end{align}
Applying \lref{L:1} for $p=2$, and with $g=\nabla {\phi}$ and $f={\phi}$, we obtain
\begin{align}\label{I:bound}
    |I|&\le C\Sob{\Lam^{\be-1}\nabla {\phi}}{L^\infty}\Sob{\phi}{L^2}\Sob{{\phi}-\tht}{L^2}\le C\left(\frac{\veps}{\veps+s}\right)^{1/2}\Sob{\phi-\tht}{L^2},
\end{align}
where we applied \cref{P:1} for $\nabla \Lam^{\be-1}\phi$ and used the fact that $\Sob{\phi}{L^2}\le C\Sob{\phi_0}{L^2}\le C$. Similarly, applying \lref{L:1} for $p=2$, and with $g=\nabla {\phi}$ and $f={\phi}-\tht$, using \cref{P:1}, and \eqref{energy:inequality}, we obtain
\begin{align}\label{J:bound}
    |J|&\le C\Sob{\Lam^{\be-1}\nabla {\phi}}{L^\infty}\Sob{{\phi}-\tht}{L^2}^2\le C\left(\frac{\veps}{\veps+s}\right)^{1/2}\Sob{\phi-\tht}{L^2}.
\end{align}
From the bounds in \eqref{dissipative:term:bound}, \eqref{I:bound}, and \eqref{J:bound}, we have that for every fixed $t$
\begin{align*}
    \Sob{\tht(t)-\phi(t)}{L^2}^2\le \Sob{\tht(0)-\phi(0)}{L^2}^2+C\int_0^t\left(\nu+\left(\frac{\veps}{\veps+s}\right)^{1/2}\right)\Sob{\tht(s)-\phi(s)}{L^2}\,ds.
\end{align*}
Applying the nonlinear Gronwall inequality as given in \cite{MitrinovicPecaricFinkBook1991} (p. 360) with $a(s)=0$, $b(s)=C\left(\nu+\left(\frac{\veps}{\veps+s}\right)^{1/2}\right)$, $p=1/2$, $u(s)=\Sob{\tht(s)-\phi(s)}{L^2}^2$, we obtain
\begin{align*}
    \Sob{\tht(t)-\phi(t)}{L^2}^2\le \left(\Sob{\tht(0)-\phi(0)}{L^2}+\frac{1}{2}C(\nu t+2\sqrt{\veps}\sqrt{\veps+t})\right)^2.
\end{align*}
Letting $\nu, \veps$ tend to zero yields the claimed result.\qed
\subsection{Proof of \cref{T:2}}
For convenience, we drop the superscript $\nu$. The first part of the proof is analogous to the proof of \cref{T:1}. In particular, the proofs of \cref{P:2} and the initial part of \cref{T:1} (up to inequality \eqref{dissipative:term:bound}) remain the same, but now applied to $\phi=\bar{\tht}$. To estimate the remaining terms of \eqref{Remainder}, we proceed as follows
\begin{align}
    &\int_0^t\int_{\mathbb{R}^2}(\bdy_s \bar{\tht}+u\cdot \nabla \bar{\tht})(\bar{\tht}-\tht)+\frac{1}{\veps}(\mathcal{R}_1 \tht)(\bar{\tht}-\tht)\,dx\,ds\notag\\&=\int_0^t\int_{\mathbb{R}^2}(u-\bar{u})\cdot \nabla \bar{\tht}(\bar{\tht}-\tht)+\frac{1}{\veps}(\mathcal{R}_1(\tht-\bar{\tht}))(\bar{\tht}-\tht)\,dx\,ds\notag\\
    &=\int_0^t\int_{\mathbb{R}^2}(u-\bar{u})\cdot \nabla \bar{\tht}(\bar{\tht}-\tht)\,dx\,ds
\end{align}
by \cref{R:3}. Proceeding just like in \eqref{eq:commutator:decom:2}, we obtain
\begin{align}
    \int_{\mathbb{R}^2} (\bar{u}-u)\cdot \nabla \bar{\tht}(\bar{\tht}-\tht)\,dx=-\frac{1}{2}\int_{\mathbb{R}^2}\left([\nabla^{\perp}\Lam^{\be-2}\cdot,\nabla \bar{\tht}](\bar{\tht}-\tht)\right)(\bar{\tht}-\tht)\,dx.
\end{align}
Applying \lref{L:1} for $p=2$, and with $g=\nabla \bar{\tht}$ and $f=\bar{\tht}-\tht$, we obtain
\begin{align}\label{nonlinear:bound:2}
    \bigg{|} \int_{\mathbb{R}^2} (\bar{u}-u)\cdot \nabla \bar{\tht}(\bar{\tht}-\tht)\,dx\bigg{|}&\le C\Sob{\Lam^{\be-1}\nabla \bar{\tht}}{L^\infty}\Sob{\bar{\tht}-\tht}{L^2}^2\notag\\& \le C\Sob{\bar{\tht}-\tht}{L^2}^2,
\end{align}
where we used the fact that $\bar{\tht}\in L^\infty(0,T^*;H^m)$ for $m>1+\be$.
The remaining term is treated just like in \cite{KosloffNichePlanas2019}. Observe that for $m>1+\be$ and $0<\al<1$, it follows that $H^m \subset H^{2\al}$.  Using \eqref{dissipative:term:bound} with $\phi=\bar{\tht}$ and applying the Young's inequality, we have
\begin{align}\label{dissipative:term:bound:2}
    \int_0^t \int_{\mathbb{R}^2} \nu \Lam^\al \bar{\tht} (\Lam^\al \bar{\tht}-\Lam^\al \tht)\,dx\,ds&\le \nu \int_0^t \Sob{\Lam^{2\al}\bar{\tht}}{L^2}\Sob{\bar{\tht}-\tht}{L^2}\,ds\notag\\&\le \frac{\nu}{2}\left(C+\int_0^t \Sob{\bar{\tht}-\tht}{L^2}^2\,ds\right)
\end{align}
for some large positive constant $C$.
\par
Using the bounds in \eqref{nonlinear:bound:2} and \eqref{dissipative:term:bound:2} in \eqref{Remainder}, we obtain
\begin{align*}
    \Sob{\bar{\tht}(t)-\tht(t)}{L^2}^2\le C\nu+\Sob{\bar{\tht}(0)-\tht(0)}{L^2}^2+(C+\nu)\int_0^t \Sob{\bar{\tht}(s)-\tht(s)}{L^2}^2\,ds.
\end{align*}
Applying the (linear) Gronwall inequality, we have
\begin{align*}
    \Sob{\bar{\tht}(t)-\tht(t)}{L^2}^2\le C(\nu+\Sob{\bar{\tht}(0)-\tht(0)}{L^2}^2)e^{((C+\nu)t)}.
\end{align*}
Letting $\nu$ tend to zero gives us the desired result.\qed
\subsection{Proof of \cref{T:3}}
To simplify notations, we drop the superscript $\nu$. In the absence of dispersive forcing, following the same steps as in \cref{P:2}, we obtain the following relative energy inequality, for all $t\in [0,T]$ and for all $\phi \in C^\infty_c ([0,T]\times \mathbb{R}^2)$
\begin{align*}
    &\mathcal{E}(\tht,\phi)\bigg{|}_0^t+\nu \int_0^t\int_{\mathbb{R}^2}|\Lam^\al (\tht-\phi)|^2\,dx\,ds\notag\\&\le \int_0^t \int_{\mathbb{R}^2}(\bdy_s \phi+u\cdot \nabla \phi)(\phi-\tht)+\nu\Lam^\al \phi (\Lam^\al \phi-\Lam^\al \tht)\,dx\,ds.
\end{align*}
We apply this inequality to the case $\phi=\tilde{\tht}$ and estimate the right-hand side just like in the proof of \cref{T:2} to obtain the desired result.\qed
\appendix

\renewcommand{\thesection}{\null}
\renewcommand{\theequation}{A.\arabic{equation}}
\section{Proof of \cref{P:2}}\label{sec:app}
We have
\begin{align*}
&\int_{\mathbb{R}^2}
\frac12|\theta-\phi|^2\,dx\Big|_0^t
+
\nu
\int_0^t
\int_{\mathbb{R}^2}
|\Lambda^\alpha(\theta-\phi)|^2
\,dx\,ds\\
&={}
\int_{\mathbb{R}^2}
\frac12(\theta^2+\phi^2)\,dx\Big|_0^t
-
\int_{\mathbb{R}^2}
\theta\phi\,dx\Big|_0^t
+
\nu
\int_0^t
\int_{\mathbb{R}^2}
\left(
(\Lambda^\alpha\theta)^2
+
(\Lambda^\alpha\phi)^2
-
2\Lambda^\alpha\phi
\Lambda^\alpha\theta
\right)
dx\,ds
\\
&={}
\int_{\mathbb{R}^2}
\frac12(\theta^2+\phi^2)\,dx\Big|_0^t
-
\int_0^t
\int_{\mathbb{R}^2}
(\partial_s\phi+u\cdot\nabla\phi)\theta
-
\nu\Lambda^\alpha\phi\Lambda^\alpha\theta
-
\frac1{\varepsilon}(\mathcal{R}_1\theta)\phi
\,dx\,ds
\\&
+
\nu
\int_0^t
\int_{\mathbb{R}^2}
\left(
(\Lambda^\alpha\theta)^2
+
(\Lambda^\alpha\phi)^2
-
2\Lambda^\alpha\phi
\Lambda^\alpha\theta
\right)
dx\,ds
\\
&={}
\int_{\mathbb{R}^2}
\frac12(\theta^2+\phi^2)\,dx\Big|_0^t
-
\int_0^t
\int_{\mathbb{R}^2}
(\partial_s\phi+u\cdot\nabla\phi)\theta
+
\nu\Lambda^\alpha\phi\Lambda^\alpha\theta
-
\frac1{\varepsilon}(\mathcal{R}_1\theta)\phi
\,dx\,ds
\\
&
+
\nu
\int_0^t
\int_{\mathbb{R}^2}
\left(
(\Lambda^\alpha\theta)^2
+
(\Lambda^\alpha\phi)^2
\right)
dx\,ds .
\end{align*}

Since, by $\operatorname{div}u=0$, we have
\[
\int_{\mathbb{R}^2}(u\cdot\nabla\phi)\phi\,dx=0,
\]
we can write
\begin{equation*}
\int_{\mathbb{R}^2}
\frac12\phi^2\,dx\Big|_0^t
=
\int_0^t
\int_{\mathbb{R}^2}
(\partial_s\phi+u\cdot\nabla\phi)\phi
\,dx\,ds,
\end{equation*}

Using the fact that the dispersive term does not contribute to the energy decay(see \cref{R:3}) along with \eqref{energy:inequality}, we obtain

\begin{equation}\label{A1}
\int_{\mathbb{R}^2}
\frac12(\theta^2+\phi^2)\,dx\Big|_0^t
+
\nu
\int_0^t
\int_{\mathbb{R}^2}
(\Lambda^\alpha\theta)^2\,dx\,ds
\le
\int_0^t
\int_{\mathbb{R}^2}
\left(
(\partial_s\phi+u\cdot\nabla\phi)\phi
-
\frac1{\varepsilon}(\mathcal{R}_1\theta)\phi
\right)
dx\,ds.
\end{equation}

If we now add to both sides of \eqref{A1} the terms

\begin{equation*}
\int_0^t
\int_{\mathbb{R}^2}
\left(
-(\partial_s\phi)\theta
-(u\cdot\nabla\phi)\theta
-\nu\Lambda^\alpha\phi\,\Lambda^\alpha\theta
+\frac1{\varepsilon}(\mathcal{R}_1\theta)\phi
+\nu(\Lambda^\alpha\phi)^2
\right)
dx\,ds,
\end{equation*}

the resulting expression gives us the desired inequality. \qed


\bibliographystyle{plain}
\bibliography{main_bib.bib}
\vspace{.3in}
\end{document}